\documentclass[11pt]{article}
\usepackage{amsfonts, amsmath, amsthm, amssymb}
\usepackage[all]{xy}

\newtheorem{theorem}{Theorem}[section]
\newtheorem{lemma}[theorem]{Lemma}
\newtheorem{cor}[theorem]{Corollary}

\newtheorem{definitiontemp}[theorem]{Definition}
\newenvironment{definition}{\begin{definitiontemp}
\normalfont}{\end{definitiontemp}}

\theoremstyle{remark}
\newtheorem*{remark}{Remark}

\newcommand{\A}{\mathbb{A}}
\newcommand{\Hom}{\text{Hom}}
\newcommand{\lb}{\Lambda_B}
\newcommand{\zp}{\mathbb{Z}_{\!p}}
\newcommand{\fp}{\mathbb{F}_{\!p}}
\newcommand{\dind}[2]
{\genfrac{}{}{0pt}{}{#1}{#2}}

\begin{document}
\title{$p$-adic superspaces and Frobenius.}
\author{A. Schwarz\thanks{Partly supported by NSF grant No. DMS 0505735.}
, I. Shapiro}
\date{}
\maketitle

\begin{abstract}
The notion of a $p$-adic superspace is introduced and used to give a
transparent construction of the Frobenius map on $p$-adic cohomology
of a smooth projective variety over $\zp$ (the ring of $p$-adic
integers), as well as an alternative construction of the crystalline
cohomology of a smooth projective variety over $\fp$ (finite field
with $p$ elements).

\end{abstract}
\section{Introduction.}\label{DPandFrob}
If $X$ is a smooth projective variety over $\mathbb{Z}$ or, more
generally, over the ring of $p$-adic integers $\zp$ one can define
the Frobenius map on the de Rham cohomology of $X$ with
coefficients in $\zp$ \cite{bando}.  This map plays an important
role in arithmetic geometry (in particular it was used in the
Wiles' proof of Fermat's Last Theorem); recently it was used to
obtain interesting results in physics \cite{inst,instanton}.
However, the construction of this map is not simple, the usual
most invariant approach is based on the consideration of the
crystalline site \cite{bando}. In any variation one uses the
notion of a DP-ideal, that is an ideal $I$ in a ring $A$ with the
key property that for $x\in I$, $x^n/n!$ makes sense.  To be
precise, one assumes the existence of operations $\gamma_n :
I\rightarrow A$, for $n\geq 0$, that mimic the operations
$x\mapsto x^n/n!$ and satisfy the same conditions (for instance
$n!\gamma_n(x)=x^n$). The ring $A$ is then called a DP-ring (DP
stands for divided powers), and a DP-morphism is a ring
homomorphism compatible with the DP-structure.

The advantage of the crystalline cohomology (\cite{illusie} is a
good review) of a scheme $X$ over $\fp=\zp/p\zp$ is that the
coefficients of the theory are in $\zp$, though the original $X$ was
defined over $\fp$. Furthermore the action of the Frobenius
endomorphism exists in this theory. DP neighborhoods play an
essential role  in defining crystalline cohomology; roughly
speaking a DP neighborhood $\widetilde{X}$ of $X$ in $Y$ is
described locally by a pair $(\hat{B},\hat{I})$ where $\hat{B}$ is
the ring of functions on $\widetilde{X}$ and $\hat{I}$ is the ideal
of the subvariety $X$, the important requirement is that $\hat{I}$
is in fact a DP ideal. However, to construct crystalline cohomology
and to analyze its relation to de Rham cohomology some technical
problems must be overcome. We will show that using the ideas of
supergeometry\footnote{A reader who is unfamiliar with
supermathematics is encouraged to glance through Sec.
\ref{supermath}.} we can make the exposition less technical, but
still completely rigorous. (Idea: Grassmann rings have divided
powers naturally and at the same time there are enough of them to
``feel" the entire DP neighborhood.) Thus in an appropriate
(Grassmann) setting the standard notion of infinitesimal
neighborhood replaces the DP envelope.

To summarize, the goals of this paper are as follows.  Use
supergeometry to give an alternative definition of crystalline
cohomology of a smooth projective variety over $\fp$.  And secondly,
to define a lifting of the action of Frobenius to the usual de Rham
cohomology of a smooth projective variety over $\zp$. Our
construction of the Grassmann neighborhood may also be applied in
the case when $X$ is not smooth over $\fp$, but in this case it is
not likely that our cohomology coincides with the crystalline one.
However, the crystalline cohomology is known not to give good
answers in the non-smooth (over $\fp$) case anyway, see \cite{ogus}
for a better cohomology.  By playing with our definitions it may be
possible to create a theory that gives good results in the
non-smooth case.

Our considerations are based on the notion of a $p$-adic
superspace defined as a covariant functor on an appropriate
subcategory of the category of $\zp$-algebras. It seems that this
notion is interesting in itself; one can hope that it can be used
to introduce and analyze ``$p$-adic supersymmetry" and ``$p$-adic
superstring" making contact with the $p$-adic $B$-model of
\cite{inst} and $p$-adic string theory (see \cite{pstrings} for a
review).

The notion of a $p$-adic superspace that we use is very general;
it is impossible to obtain any significant results in such
generality.  In all of our examples however, we consider functors
defining a superspace that are prorepresentable in some sense.  We
sketch the proof of this fact in the Appendix and show how to use
it  to get a more conceptual derivation of some of our
statements.

Let us note that a lot of what follows does not work for the even
prime and so we omit that case by default.  Some modifications
designed to allow for the even prime are mentioned in Sec.
\ref{evenprime}.

Finally,  any functor, by default,  is a covariant functor.  All the
rings are assumed to be unital.  All varieties are of finite type
over $\zp$.

\subsection{Summary of main definitions and results.}

The category $\Lambda$ that serves as the source for most of our
functors is described in Definition \ref{maindef}.  The main notion
in this paper is that of a $p$-adic superspace (see Definition
\ref{padicsuperspace}).  The replacement for a DP neighborhood in
$\mathbb{P}^n_{\zp}$ of $W\subset \mathbb{P}^n_{\fp}$, namely
$\widetilde{W}$ (called the Grassmann neighborhood), is given in
Definition \ref{dpnhood}. Finally, our notion of de Rham cohomology
of a $p$-adic superspace is summarized in Definition
\ref{derhamcoh}.

The appendix contains a discussion of prorepresentable $p$-adic
superspaces and should be considered as the correct general setting
for this paper.

For the Frobenius action on $\widetilde{W}$ please see Lemma
\ref{frobenius}.  The investigation of functions on $\widetilde{pt}$
(where $pt$ is a point in a line or $\mathbb{P}^1$) is undertaken in
Sec. \ref{pita} and summarized in Corollary \ref{pitaproof}.

The key points of the paper are the comparison results between our
de Rham cohomology of $p$-adic superspaces and the usual de Rham
cohomology as well as the crystalline cohomology.  Namely, for a
smooth $V$ over $\zp$, Lemma \ref{bridge} relates the de Rham
cohomology of the completion of $V$ with respect to $p$ and the de
Rham cohomology of the $p$-adic superspace associated to $V$.  Lemma
\ref{derhamtodr} shows that when $V$ is also projective, the
completion can be removed.  The main theorem of the paper is Theorem
\ref{main} that establishes the isomorphism between the de Rham
cohomology of a smooth projective $V$ over $\zp$ and the de Rham
cohomology of the Grassmann neighborhood in $\mathbb{P}^n$ of the
$p$-adic superspace associated to $V$.  Theorem \ref{crysderham}
compares the crystalline cohomology of $W$, a smooth projective variety over
$\fp$, with the de Rham cohomology of its Grassmann neighborhood (in $\mathbb{P}^n_{\zp}$)
$\widetilde{W}$.

Theorem \ref{main} is refined in Theorem \ref{mainprime} where it is
shown that the isomorphism is in fact compatible with the Hodge
filtration on one side and a very natural filtration on the other.
This allows us to re-prove some divisibility estimates for the
action of Frobenius on the de Rham cohomology of $V$.  The most
useful of these is stated in Lemma \ref{divest}.

The main motivation for this paper is Corollary \ref{frobactoncoh}
that shows that a smooth projective variety over $\zp$ has an action
of Frobenius on its de Rham cohomology.

\subsection{Supermathematics.}\label{supermath}
The present paper aims to explain crystalline cohomology from the
point of view of supergeometry.  However the  knowledge of
supergeometry itself is not essential  for the reading of this text
(at least formally).  Our treatment of spaces as functors allows for
a quick jump from the familiar commutative setting to the
supercommutative one.

The functorial approach to the definition of the superspace (in an
essentially different form) was advocated in \cite{ssup} and
\cite{konech}. It is analogous to the  functor approach to the
theory of schemes in algebraic geometry.  An excellent reference for
this point of view is \cite{jantzen} where this theory is developed
from the very beginning; groups and Lie algebras are treated as
well. This treatment of algebraic geometry generalizes immediately
to the
 setting of supergeometry by replacing commutative rings with
supercommutative ones.  The reader is invited to consult these
papers in the case of necessity. For a more complete understanding
of supermathematics we recommend \cite{super1}, as well as
\cite{super2} for those looking for a more abstract and conceptual
picture.

The very few concepts that we do need are explained  below.

\begin{definition}
By a supercommutative ring $A$ we mean a
$\mathbb{Z}/2\mathbb{Z}$-graded ring (i.e. $A=A^0\oplus A^1$ with
multiplication respecting the grading) such that $$a_i\cdot
a_j=(-1)^{ij}a_j\cdot a_i$$ where $a_i\in A^i$.\footnote{This is the
sign rule. Whenever two odd objects are to be written in the
opposite order, a minus sign is the price.  This is the mildest
possible modification of the usual commutativity.} We say that $a\in
A^0$ is \emph{even} and $b\in A^1$ is \emph{odd}.
\end{definition}

\begin{remark}
Any commutative ring $A$ is also supercommutative with $A^1=0$.  We
often use freely generated rings, i.e. $B[x_i,\xi_j]$ with $B$
commutative, and generators $x_i$ commuting ($x_i x_j=x_j x_i$) and
$\xi_j$ anti-commuting ($\xi_i\xi_j=-\xi_j\xi_i$), that is $x_i$ are
even and $\xi_i$ are odd, thus $x_i\xi_j=\xi_j x_i$.  Note that
specifying the parity of $B$, $x_i$ and $\xi_i$ is sufficient to
define the $\mathbb{Z}/2\mathbb{Z}$-grading. It
is convenient to consider $x_i$ and $\xi _j$ as commuting and anticommuting variables.
\end{remark}

\begin{definition}
Denote by $Super$ the category of supercommutative rings, with
morphisms being algebra homomorphisms respecting the grading (i.e. a
morphism preserves the $\mathbb{Z}/2\mathbb{Z}$ decomposition). In
other words, morphisms are parity preserving homomorphisms.
\end{definition}

\begin{remark}
To give the reader a feeling for the morphisms in $Super$ we point
out that $\Hom_{Super}(\mathbb{Z}[x_1,...,x_n], A)=(A^0)^n$,
$\Hom_{Super}(\mathbb{Z}[\xi_1,...,\xi_n], A)=(A^1)^n$ and
$\Hom_{Super}(\mathbb{Z}[x_1,...,x_n,\xi_1,...,\xi_m],
A)=(A^0)^n\times(A^1)^m$.  We note the standard use of $x_i$ for
even variables and $\xi_i$ for odd.
\end{remark}

The advantage of the functorial approach to spaces is that it is
deceptively simple.  We can say that a superspace, in the most
general sense, is a functor from $Super$ (or from a subcategory of
$Super$) to $Sets$ (the category of sets). Morphisms are natural
transformations of functors.  This definition ignores the most
important aspect of spaces, namely that they should be local. This
can be included in the definition in many ways. For the general
method see \cite{jantzen} and compare with our  approach based on
the consideration of the body of a superspace (Definition
\ref{bodyof}) and prorepresentability (Appendix). Supergroups are
functors from $Super$ to $Groups$ (the category of groups), etc. The
action of a supergroup $G$ on a superspace $X$ can be defined very
naturally: for a supercommutative ring $A$  we  should have an
action of the group $G(A)$ on the set $X(A)$ satisfying some
conditions of functoriality.

We define the rest of the concepts
as needed.

In fact the reader is already quite familiar with some of the
supermathematics that we use.  For example, the de Rham complex of a
smooth affine variety $V$ is a supercommutative algebra with a
differential.  This algebra can be interpreted as the algebra of
functions on the odd tangent space $\Pi TV$ of $V$ and the
differential as an action of the odd line $\mathbb{A}^{0,1}$ on $\Pi
TV$.  This is thoroughly discussed in Sec. \ref{De Rham}.

\subsection{Main constructions.}

Let us explain our constructions in some detail. (The reader can
skip this explanation. However it could be useful for some people,
in particular, for readers with a background in mathematical
physics.)

If a (projective or affine) variety is specified by means of
equations with coefficients belonging to a commutative ring $R$ it makes sense
to consider the unknowns as belonging to any $R$-algebra. This means
that every variety of this kind (variety over $R$) specifies a
functor on the category of $R$-algebras with values in the category
of sets.

This remark prompts a preliminary definition of a superspace over a
ring $R$ as a functor taking values in the category of sets and
defined on the category of $\mathbb{Z}/2\mathbb{Z}$-graded
supercommutative $R$-algebras or better yet on its subcategory
$\Lambda$. (Morphisms are parity preserving homomorphisms.) We
define a map of superspaces as a map of functors; this definition
depends on the choice of the category $\Lambda$. In particular, in
the case when $R$  is a ring of $p$-adic integers $\zp$ we take as
$\Lambda$ the category of supercommutative rings of the form
$B\otimes \Lambda_n$, where $B$ is a  finitely generated commutative
ring in which the multiplication by $p$ is nilpotent and $B/pB$ does
not contain nilpotent elements. Here $\Lambda_n$ stands for the
Grassmann ring (free supercommutative ring with  $n$ odd
generators.) Every ring in the category $\Lambda$ can be considered
as a $\zp$-algebra since multiplication by a series $\sum_n a_n p^n$
is well defined because multiplication by $p$ is nilpotent.
Considering functors on the above category we come to the definition
of a $p$-adic superspace.

The typical example of a  superspace  is ${\A}^{n|k,m}$
corresponding to a functor, that sends an $R$-algebra to a set of
rows with $n$ even, $k$ even nilpotent elements and $m$ odd elements
of the algebra. A function on a superspace is a map of this
superspace to ${\A}^{1,1}={\A}^{1|0,1}$. In the $p$-adic case the
functions on ${\A}^{n|k,m}$ are described in Theorem
\ref{theorem:ankm}. In particular, the functions on  ${\A}^{0|1,0}$
correspond to series with divided powers $\sum a_n{x^n\over n!}$
where $a_n\in\zp$. (Notice that in this statement it is important
that the functors are defined on the category $\Lambda$ described
above.  We could consider the functors on the larger category of all
supercommutative $\zp$-algebras;  then only a series of the form
$\sum a_n x^n$ with $a_n\in \zp$ corresponds to a function.)

The definition of a superspace in terms of a functor is too general;
one should impose some additional restrictions to develop an
interesting theory. One of the possible ways is to impose some
conditions on the ``bosonic part" of the superspace (i.e., on the
restriction of the functor to commutative $R$-algebras.) In
particular, we notice that the category $\Lambda$ we used in the
$p$-adic case contains the category of commutative $\fp$-algebras
without nilpotent elements; in this case the restriction of the
functor specifying a $p$-adic superspace $Y$ should correspond to a
variety over $\fp$ (the body $[Y]$ of the superspace $Y$).

For any $C\in\Lambda$, let us set $\rho(C)=C/C^{nilp}$, where
$C^{nilp}$ is the ideal of nilpotent elements in $C$. In our case
$C^{nilp}$ is generated by $p$ and the odd generators.  From
definitions we see that $\rho(C)$ is an $\fp$-algebra.  Thus the
natural projection $C\rightarrow \rho(C)$ induces a map
$\pi:Y(C)\rightarrow [Y](\rho(C))$.  For every (open or closed)
subvariety $Z\subset[Y]$ we can define a $p$-adic subsuperspace
$Y|_Z\subset Y$ as a maximal  subsuperspace of $Y$ having $Z$ as its
body. More explicitly, $Y|_Z(C)=\pi^{-1}(Z(\rho(C)))$; we can say
that $Y|_Z$ is a  subsuperspace  of $Y$ over $Z$.

We give a definition of the body only in $p$-adic case, but similar
constructions also work in other situations.

For every superspace $Y$ over a ring $R$ we can introduce a notion
of a differential form and of the de Rham differential.
(Differential forms are defined as functions on the  superspace $\Pi
TY$ that parameterizes the maps from the superspace
$\A^{0,1}=\A^{0|0,1}$ to $Y$.) One can try to define the cohomology
of $Y$ as cohomology of the differential $R$-module $\Omega(Y)$ of
differential forms on $Y$, but this definition does not capture the
whole picture (it corresponds to the consideration of the Hodge
cohomology $H^{k,0}$).  The right definition of de Rham cohomology
of $Y$ can be given in terms of hypercohomology of a sheaf of
differential $R$-modules on the body of $Y$. (To an open subset
$Z\subset[Y]$ we assign the module $\Omega(Y|_Z)$ of differential
forms on $Y|_Z$.)

Now we are able to define our analogue of the crystalline cohomology
of a projective $\fp$-variety $X\subset \mathbb{P}^n_{\fp}$ where
$\mathbb{P}^n_{\fp}$ stands for projective space over $\fp$. One can
regard $\mathbb{P}^n_{\fp}$ as a body of the $p$-adic projective
superspace $\mathbb{P}^n$; this remark permits us to consider
$\mathbb{P}^n|_X$ (the subsuperspace of $\mathbb{P}^n$ over $X$). We
define ``crystalline" cohomology of $X$ as the de Rham  cohomology
of the $p$-adic superspace $\mathbb{P}^n|_X$. The Frobenius map $Fr$
acts naturally on this  cohomology. (The usual action of $Fr$ on
$\mathbb{P}^n$ sending every homogeneous coordinate to its $p$-th
power preserves the $\fp$-variety $X$ and therefore
$\mathbb{P}^n|_X$.) We prove that for a smooth $\fp$-variety $X$ the
cohomology of $\mathbb{P}^n|_X$ coincides with the $p$-adic de Rham
cohomology of any variety $X'$ over $\zp$ that gives $X$ after
reduction to $\fp$ (and therefore the Frobenius acts on the
cohomology of $X'$); see Corollary \ref{frobactoncoh}.

\section{Category $\Lambda$.}

Consider the local ring $\zp$ with the maximal ideal $p\zp$ (it is a
DP ideal since $p^n/n!$ which is obviously in $\mathbb{Q}_p$ is
actually in $\zp$ because $\text{ord}_p n!\leq n$). Our main object
is the category $\Lambda$.

\begin{definition}\label{maindef}
Let $\Lambda$ be the category with objects $\lb$ that are super
commutative rings freely and finitely generated over a commutative
ring $B$\footnote{A typical example of such $B$ is
$\mathbb{Z}/p^n\mathbb{Z}$.} ($B$ is allowed to vary) by odd
generators. More precisely, we require that $B$ is a finitely
generated commutative ring such that $p$ is nilpotent, and $B/pB$
has no nilpotent elements. We write the $\mathbb{Z}/2\mathbb{Z}$
grading of $\lb$ as follows $\lb=\lb^{0}\oplus\lb^{1}$ where
$\lb^{0}$ is even and $\lb^{1}$ is odd. The morphisms are parity
preserving homomorphisms.
\end{definition}

The requirement that $p$ be nilpotent is necessary to allow for
evaluation of infinite series such as $\sum a_n x^n$ at elements
$pb$ with $b\in B$. Also, for example, we see that $\lb$ is a
$\zp$-module (in fact it is a $\mathbb{Z}/p^N\mathbb{Z}$-module for
an $N$ large enough). This means that we can consider $\Lambda$ as a subcategory of the category of $\zp$-algebras.

\begin{remark}
We use the notation $\lb$ to denote a generic element of $\Lambda$
to emphasize the fact that it is $B$ that is important, not the
number of odd variables.  The notation is meant to remind the reader
of the exterior algebra (Grassmann algebra) where the field has been
replaced with $B$. A more precise statement would be that $\lb\simeq
B[\xi_i]$ where $B[\xi_i]$ is a polynomial ring in $\xi_i$s with
coefficients in $B$, however in this case the variables $\xi_i$ are
anti-commuting, that is $\xi_i\xi_j=-\xi_j\xi_i$.  The parity is
determined by the total number of variables $\xi_i$.

\end{remark}

Denote by $\lb^+$ the ideal in $\lb$ generated by $\xi_1,...,\xi_n$.
Notice that $\Lambda$ contains the category of $\fp$-algebras
without nilpotent elements as a full subcategory, and there is a
retraction onto it that sends $\lb$ to $\lb/(pB+\lb^+)=B/pB$.  The
ideal $pB+\lb^+$ is to play a very important role for us. One should
mention that it can be characterized by the fact that it consists
exactly of the nilpotent elements of $\lb$ (this follows from the
lack of nilpotent elements in $B/pB$).  As a consequence we see that
if $A\rightarrow \lb$ is any morphism and $I\subset A$ is a
nilpotent ideal, then $I$ is carried to $pB+\lb^+$ by the morphism.
However its most important property is explained in the following
Theorem.

\begin{theorem}\label{dptheorem}
For every $\lb\in\Lambda$, $pB+\lb^+\subset\lb$ is naturally a DP
ideal, i.e. there are operations $\gamma_n: pB+\lb^+\rightarrow \lb$
that satisfy the axioms in \cite{bando}.  Furthermore, any morphism
in $\Lambda$ preserves this structure automatically.
\end{theorem}

\begin{remark}
Thus every object in $\Lambda$ is in fact a DP pair $(\lb,pB+\lb^+)$
and any morphism preserves the DP structure.  This explains our
choice of $\Lambda$.
\end{remark}

\begin{proof}
Represent $B$ as $\zp[x_i]/I$, where $\zp[x_i]$ is the polynomial
algebra over $\zp$. Note that since $\zp$ is torsion free
$\Lambda_{\zp[x_i]}\subset\Lambda_{\mathbb{Q}_p[x_i]}$ and we can
define
$\gamma_n:\Lambda_{\zp[x_i]}\rightarrow\Lambda_{\mathbb{Q}_p[x_i]}$
by $\gamma_n(x)=x^n/n!$ for all $n$.  We claim that
$p\zp[x_i]+\Lambda_{\zp[x_i]}^+$ maps under $\gamma_n$ to
$\Lambda_{\zp[x_i]}$. (Thus $\gamma_n$ define a DP structure on the
pair $(\Lambda_{\zp[x_i]},p\zp[x_i]+\Lambda_{\zp[x_i]}^+)$.) To show
this it is sufficient to check that $x^n/n!$ is in
$\Lambda_{\zp[x_i]}$ for $x\in p\zp[x_i]$ and for $x\in
\Lambda_{\zp[x_i]}^+$.

For $x\in p\zp[x_i]$, this follows from the observation that
$p^n/n!\in\zp$. Now suppose that $e\in \Lambda_{\zp[x_i]}^{0+}$, let
$e=e_1+...+e_k$ where $e_i$ are even and of the form
$f_i\xi_{i_1}...\xi_{i_s}$, i.e. write $e$ as the sum of monomials
in $\xi_i$s.  Notice that $e_i^n=0$ if $n>1$. So that
$$e^n=(e_1+...+e_k)^n=\sum_{\sum n_i=n}\dfrac{n!}
{n_1!...n_k!}e_1^{n_1}...e_k^{n_k} =\sum_{\dind{\sum
n_i=n}{n_i=0\,\text{or}\,1}}n!e_1^{n_1}...e_k^{n_k}.$$   Consider an
element $x=e+o\in \Lambda_{\zp[x_i]}^+$  with $e$ even and $o$ odd.
Then $x^n=(e+o)^n=e^n+n e^{n-1} o$ and we are done.

Since $\gamma_n$ satisfy the axioms for a DP structure we obtain in
this way a DP structure on the pair
$(\Lambda_{\zp[x_i]},p\zp[x_i]+\Lambda_{\zp[x_i]}^+)$.  Note that
$\lb=\Lambda_{\zp[x_i]}/I[\xi_j]$, thus it inherits a DP structure
from $\Lambda_{\zp[x_i]}$ if (and only if) $I[\xi_j]\cap
(p\zp[x_i]+\Lambda_{\zp[x_i]}^+)$ is preserved by $\gamma_n$.  But
 $I[\xi_j]\cap (p\zp[x_i]+\Lambda_{\zp[x_i]}^+)=I\cap
p\zp[x_i]+\Lambda_{I}^+$, and clearly $\Lambda_{I}^+$ is preserved
by $\gamma_n$.  As for $I\cap p\zp[x_i]$, if $pf\in I$ then
$(pf)^n/n!=pfg$ with $g\in\zp[x_i]$ thus $(pf)^n/n!\in I$.  We
conclude that $\lb$ inherits a DP structure on $pB+\lb^+$.

Next we observe that the DP structure obtained as above does not
depend on a particular representation of $B$ as a quotient of a
polynomial algebra. Namely, if we represent $B$ as $\zp[y_j]/J$ and
obtain a DP structure on $\lb$ in that way, then the identity map on
$\lb$ lifts to a homomorphism from $\Lambda_{\zp[x_i]}$ to
$\Lambda_{\zp[y_j]}$ (because $\Lambda_{\zp[x_i]}$ is free) that is
automatically compatible with DP structure (since $\gamma_n$ is just
$x^n/n!$). Because the projections to $\lb$ are DP compatible by
definition, the identity map is DP compatible as well.

If $\lb\rightarrow\Lambda_C$ is any morphism then it lifts to a
morphism of the free algebras that cover $\lb$ and $\Lambda_C$ as
above. The lifting is again automatically compatible with DP
structure, ensuring that $\lb\rightarrow\Lambda_C$ is a DP morphism.

\end{proof}

Because of the nature of our definition of DP structure on $\lb$ we
use the more suggestive $x^n/n!$ instead of the more accurate
$\gamma_n$ to denote the DP operations.  As we have shown above the
symbol $x^n/n!$ is functorially defined for elements of the ideals
$pB+\lb^+$.

\begin{remark}
Note that a sufficient condition on $J$ for $\Lambda_{\zp[x_i]}/J$
to inherit a DP structure is that it be a $\xi$-homogeneous ideal,
i.e. $J=\bigoplus_\alpha J\cap\zp[x_i]\xi^\alpha$, where $\alpha$ is
a multi-index.  An example used in the Theorem above is
$J=I[\xi_j]$.  Perhaps one can use this observation to enlarge the
category $\Lambda$.
\end{remark}

\begin{remark}
Given any super-commutative $\zp$-algebra $A=A^0\oplus A^1$, we may
define $A^+\subset A$ by $A^+=pA+A^1 A$, generalizing the ideal
$pB+\lb^+\subset\lb$. This ideal is functorial, however it is not
clear why it should have any DP structure.  Various conditions may
be imposed to ensure this.  The previous remark is an example.
\end{remark}

\section{Superspaces, neighborhoods and
Frobenius.}\label{supspnhoodsandfrob}

We would like to base our definition of a $p$-adic superspace on
the notion of a functor from $\Lambda$ to $Sets$, the category of
sets. Since we are interested in studying geometric objects, we
would like to impose conditions that would make the functor
``local", the easiest way to do it is through the notion of the
body of a $p$-adic superspace.

\begin{definition}\label{padicsuperspace}
A $p$-adic superspace $X$ is a functor (covariant) from $\Lambda$ to
$Sets$, such that the restriction to the full subcategory of
$\fp$-algebras without nilpotent elements corresponds to a variety
$[X]$ over $\fp$.
\end{definition}

\begin{definition}\label{bodyof}
The body of a $p$-adic superspace $X$ is the variety $[X]$.
\end{definition}

\begin{definition}
A map $\alpha: X\rightarrow Y$ of superspaces is a natural
transformation from $X$ to $Y$.
\end{definition}

A more familiar object, the purely even superspace, is obtained by
requiring that the functor factors through $\Lambda^{0}$, the
category with objects of the form $\lb^{0}$.

We have the usual functors $\mathbb{A}^n$ and
$\mathbb{P}^n$,\footnote{These functors, and others like them below,
can actually be defined as usual superspaces, i.e. they can be
obviously extended to the category of all supercommutative rings.
Here we use their restriction to $\Lambda$, but no extra structure
of $\Lambda$ is required.} where
$$\mathbb{A}^n(\lb)=\{(r_1,...,r_n)|r_i\in \lb^{0}\}$$ and
$$\mathbb{P}^n(\lb)=\{(r_0,...,r_n)|r_i\in \lb^{0},\sum\lb
r_i=\lb\}/(\lb^{0})^\times.$$ Note that these are purely even.
More generally we can define the superspace
$$\mathbb{A}^{n,m}(\lb):=\{(r_1,...,r_n,s_1,...,s_m)|r_i\in
\lb^{0},s_i\in \lb^{1} \}.$$ A further generalization that we will
need is
$$\mathbb{A}^{n|k,m}(\lb):=\{(r_1,...,r_n,t_1,...,t_k,s_1,...,s_m)|r_i\in
\lb^{0},t_i\in pB+\lb^{0+},s_i\in \lb^{1} \},$$ in other words $r_i$
are even elements, $t_i$ are even nilpotent and $s_i$ are odd
elements.  One can also define
$$\mathbb{P}^{n,m}(\lb):=\{(r_0,...,r_n,s_1,...,s_m)|r_i\in
\lb^{0},s_i\in \lb^{1},\sum\lb r_i=\lb\}/(\lb^{0})^\times$$ but we
will not need it.

\begin{remark}
The most important cases from the above are
$$\mathbb{A}^{1|0,0}(\lb)=\lb^0$$ $$\mathbb{A}^{0|1,0}(\lb)=pB+\lb^{0+}$$
$$\mathbb{A}^{0|0,1}(\lb)=\lb^1$$ they are the main building
blocks for the theory in this paper.

\end{remark}

\begin{definition}
A function on a $p$-adic superspace $X$ is a natural
transformation from $X$ to the superline $\A^{1,1}$.
\end{definition}

Considering all natural transformation from $X$ to the superline
$\A^{1,1}$ we obtain the set of functions on $X$.  It is easily
seen to be a ring by observing that the functor $\mathbb{A}^{1,1}$
takes values in the category of supercommutative rings.

A very versatile notion that we will need is that of a restriction
of a $p$-adic superspace $Y$ to a subvariety $Z$ (it need not be
open or closed) of its body.  It is the maximal subsuperspace of
$Y$ with body $Z$.  More precisely we have the following
definition.

\begin{definition}
Let $Y$ be a $p$-adic superspace and $Z$ a subvariety of $[Y]$,
then the $p$-adic superspace $Y|_Z$ is defined to make the
following diagram cartesian.
$$\xymatrix{Y|_Z(\lb)\ar@{^{(}->}[r]\ar[d]
& Y(\lb)\ar[d]^{\pi}\\
Z(B/pB)\ar@{^{(}->}[r]^i& Y(B/pB)}$$
\end{definition}

Consider the following ``local" analogue of functions on $X$.

\begin{definition}
Define the pre-sheaf of rings $\mathcal{O}_X$ on $[X]$ by setting
$\mathcal{O}_X(U)$ to be the ring of functions on $X|_U$, for any
open $U$ in $[X]$.
\end{definition}

There is no reason to expect that the pre-sheaf $\mathcal{O}_X$ is a
sheaf in general.  Thus the usual thinking about functions in terms
of coordinates is not advised. However, it is a sheaf for all the
superspaces that we consider in this paper. If one wants a more
general setting in which $\mathcal{O}_X$ is a sheaf, one should
consider prorepresentable superspaces as defined in the Appendix.

\begin{definition}
Denote by $R\left<y_1,...,y_k\right>$ the ring whose elements are
formal expressions $\sum_{K\geq 0} a_K y^K/K!$, where $a_K\in R$ and
$y_i$s are commuting variables.\footnote{Here and below $I$, $J$,
$K$, $T$ are multi-indices and $T!=t_1 ! t_2 ! ...$.}  Note that
$K!$ need not be invertible in $R$. We call
$R\left<y_1,...,y_k\right>$ the ring of power series with divided
powers.
\end{definition}

\begin{remark}
Clearly we can add such expressions, but it is also easy to see that
we can multiply them since $\left(\sum a_i y^i/i!\right)\left(\sum
b_j y^j/j!\right)=\sum\left(\sum_{i+j=n}\frac{n!}{i!j!}a_i
b_j\right)y^n/n!$.
\end{remark}

\begin{theorem}\label{theorem:ankm}
The functions on $\mathbb{A}^{n|k,m}$ are isomorphic as a ring to
$$\left\{\sum_{I,J,T\geq 0}a_{I,J,T}x^I \xi^J y^T/T!\right\}$$ where $x_i$ and
$y_i$ are even and $\xi_i$ are odd, and $a_{I,J,T}\in\zp$ with
$a_{I,J,T}\rightarrow 0$ as $I\rightarrow \infty$.
\end{theorem}

\begin{proof}
Recall that
$\mathbb{A}^{n|k,m}(\lb)=\{(r_1,...,r_n,t_1,...,t_k,s_1,...,s_m)|r_i\in
\lb^{0},t_i\in pB+\lb^{0+},s_i\in \lb^{1} \}$, then every
$\sum_{I,J,T\geq 0}a_{I,J,T}x^I \xi^J y^T/T!$ can be evaluated at
every $(r,t,s)$ by setting  $x=r$, $\xi=s$ and $y=t$ to obtain an
element of $\lb=\mathbb{A}^{1,1}(\lb)$.

Furthermore, $\mathbb{A}^{n,m}$ is pro-represented\footnote{The
phrase $F$ is pro-represented by $C_n$ is used here  in a more
narrow sense than in the Appendix. Namely, we mean  that
$F=\varinjlim h_{C_n}$, where $h_{C_n}=\text{Hom}(C_n,-)$, and $C_n$
form an inverse system of objects in the category. A good reference
on pro-representable functors (in this sense) in (non-super)
geometry is \cite{artin}.} by
$(\zp/p^N\zp)[x_1,...,x_n,\xi_1,...,\xi_m]\in\Lambda$ (they are in
$\Lambda$ because $\fp[x_i]$ has no nilpotent elements). Here $x_i$
are even and $\xi_i$ are odd.  Thus\footnote{This is an application
of the Yoneda Lemma.} the functions are
$$\displaystyle\varprojlim_N
\,(\zp/p^N\zp)[x_1,...,x_n,\xi_1,...,\xi_m]$$ i.e. of the form
$\sum_{I,J\geq 0}a_{I,J}x^I\xi^J$ with $a_{I,J}\rightarrow 0\in\zp$
as $I\rightarrow\infty$.  The case of $\A^{0|k,0}$ is not as
trivial, the issue is that it is ``pro-represented" by
$(\zp/p^N\zp)\left<y_1,...,y_k\right>$ but these are not in
$\Lambda$ (in $\fp\left<y_i\right>$ all $y_i$ are \emph{nilpotent}).
Thus the proof of the complete Theorem is postponed until it appears
as Corollary \ref{ankm}.
\end{proof}

\begin{definition}\label{vartosup}
Given a variety $V$ over  $\zp$ (which can be viewed as a functor
from the category of commutative $\zp$-algebras to $Sets$) we can
define the associated $p$-adic superspace $X_V$ by setting
$X_V(\lb)=V(\lb^{0})$.
\end{definition}

Note that information is lost in passing from the variety to its
associated superspace.  More precisely, $V$ and its $p$-adic
completion $\hat{V}_p$ will have the same associated superspace.
(See Lemma \ref{difference} below; we return to this discussion in
Sec. \ref{De Rham}.) This is best illustrated by considering the
functions on $\A^n$. As a variety over $\zp$ its functions are by
definition $\zp[x_1,...,x_n]$, however when considered as a $p$-adic
superspace one gets the much larger ring
$$\displaystyle\varprojlim_N \,(\zp/p^N\zp)[x_1,...,x_n]$$ of
functions\footnote{It consists of series with $p$-adically vanishing
coefficients.}. Observe that the functions on the purely odd affine
space $\A^{0,m}$ do not change. The crucial point for us is the
metamorphosis that the functions on $\A^{0|1,0}$ undergo, as we pass
from considering it as a variety over $\zp$ (prorepresented by
$\zp[[x]]$) to the associated superspace; they transform from power
series to divided power series. It is this observation that
motivates the present paper.

\begin{remark}
Our use of $\mathbb{A}^{n,m}$ for a $p$-adic superspace is somewhat
misleading as that symbol is standard for a superscheme; in
particular $\mathbb{A}^n$ usually denotes (with the ground ring
being implicitly $\zp$) $\text{Spec}(\zp[x_1,...,x_n])$.  There is
no confusion however if we make explicit wether we are discussing a
variety or a $p$-adic superspace associated to it.  When we need to
make the distinction  explicit, we use $X_V$ for the $p$-adic
superspace associated to a variety $V$.
\end{remark}

\begin{lemma}\label{difference}
Let $V$ be a  variety over $\zp$, then the body of $X_V$ is the
restriction  of $V$ to $\fp$, i.e., $$[X_V]=V|_{\fp}$$ and the
functions on $X_V$ (as a sheaf on $[X_V]$) are given by the
completion of the functions on $V$ at the subvariety $V|_{\fp}$,
i.e., $$\mathcal{O}_{X_V}=\widehat{(\mathcal{O}_V)}_p$$ thus making
precise the difference between $V$ and $X_V$.
\end{lemma}

\begin{proof}
That $[X_V]=V|_{\fp}$ is immediate from the definition.  The
question of functions is local, so assume $V=\text{Spec} A$.  Then
$X_V$ is pro-represented by $\varprojlim A/p^n A$, so that
$\mathcal{O}_{X_V}=\widehat{A}_p$.
\end{proof}

One of the most important notions of this paper is that of the
infinitesimal neighborhood of one $p$-adic superspace inside
another.  It is meant to replace the DP-neighborhood.

\begin{definition}\label{nhood}
Let $X\subset Y$ be $p$-adic superspaces, define $\widetilde{X}$,
the infinitesimal neighborhood of $X$ in $Y$ by
$\widetilde{X}=Y|_{[X]}$.\footnote{Note that $\widetilde{X}$ depends
only on $Y$ and $[X]$, compare with Definition \ref{dpnhood}.}
\end{definition}

\noindent\textbf{Example.} Let $\A^n\hookrightarrow\A^{n+1}$ be an
inclusion of $p$-adic superspaces, i.e. $\A^n(\lb)=\lb^0\times
...\times\lb^0\rightarrow\lb^0\times
...\times\lb^0\times\{0\}\subset\A^{n+1}(\lb)$, then it follows
directly that $\A^{n+1}|_{[\A^n]}(\lb)=\lb^0\times
...\times\lb^0\times(pB+\lb^{0+})$, so that
$\widetilde{\A^n}=\A^{n|1}$.

Suppose that $W\subset\mathbb{P}^n_{\fp}$ is a possibly non-smooth
$\fp$-variety.  We want to define the notion replacing that of a DP
neighborhood of $W$ in $\mathbb{P}^n_{\zp}$.  We do this as follows.
The inclusion of varieties over $\fp$ gives rise to a subvariety $W$
of the body of the $p$-adic superspace $\mathbb{P}^n$. Let us use
the same notation as above, namely $\widetilde{W}$ to denote
$\mathbb{P}^n|_W$, this is the infinitesimal neighborhood of $W$
that behaves much better than $W$ itself.

\begin{definition}\label{dpnhood}
We call $\widetilde{W}$ as above, the DP neighborhood of $W$ in
$\mathbb{P}^n$.  An alternative name that we sometimes use is
Grassmann neighborhood.
\end{definition}

\begin{remark}
$\widetilde{W}$ is a $p$-adic superspace whereas $W$ was a variety.
While not the same kind of object, we can nevertheless define many
notions for $p$-adic superspaces that we have for varieties. For
example as we see in the next section, we may consider the de Rham
cohomology of a $p$-adic superspace.
\end{remark}

We have the usual action of the Frobenius map $Fr$ on the $p$-adic
superspace $\mathbb{P}^n$ via raising each homogeneous coordinate to
the $p$th power. The restriction of $Fr$ to the body of
$\mathbb{P}^n$ preserves $W$ therefore we have an action of $Fr$
also on $\widetilde{W}$.  Summarizing we get:

\begin{lemma}\label{frobenius}
Let $W\subset\mathbb{P}^n$ be an inclusion of varieties over
$\mathbb{F}_p$, then any lifting of the action of Frobenius from
$\mathbb{P}^n_{\fp}$ to $\mathbb{P}^n_{\zp}$ (i.e., a choice of
homogeneous coordinates) restricts to an action of Frobenius on
$\widetilde{W}$.
\end{lemma}

\section{De Rham cohomology of superspaces.}\label{De Rham}

Let us briefly review the notion of de Rham cohomology from the
point of view of superspaces.  This point of view lends itself most
naturally to a generalization applicable in our setting.  We begin
by introducing the notion of the odd tangent space to a $p$-adic
superspace $X$.

\begin{remark}
For a $\lb\in\Lambda$, we denote by $\lb[\xi]$ the ring $\lb$ with
an adjoined extra odd variable $\xi$.  More precisely, given a
supercommutative ring $R$, we can form $R[\xi]$ by considering
expressions of the form $a+b\xi$, with multiplication defined by
$(a+b\xi)(c+d\xi)=ac+(ad+(-1)^{\bar{c}}bc)\xi$.  Where $a,b,c,d\in
R$ and $\bar{c}$ is the parity of $c$.
\end{remark}

\begin{definition}
Let $X$ be a $p$-adic superspace, define a new $p$-adic superspace
$\Pi TX$, the odd tangent space of $X$, by $$\Pi
TX(\lb)=X(\lb[\xi]).$$
\end{definition}

Functions on $\Pi TX$ will serve as the differential forms on $X$.
We will define the differential $d$ and the grading on differential
forms in terms of an action of a supergroup on $\Pi TX$.

Note that there is a natural map $\lb[\xi]\rightarrow\lb$ that sends
$\xi$ to $0$.  This defines a map of $p$-adic superspaces
$$\pi:\Pi TX\rightarrow X$$ and a corresponding map on the bodies\footnote{When $X$ is purely even, $[\Pi TX]=[X]$.}
$$[\pi]:[\Pi TX]\rightarrow [X].$$

The superspace $\Pi TX$ carries an action\footnote{These definitions
become more transparent when one thinks of $\Pi TX$ as the
superspace parameterizing the maps from $\mathbb{A}^{0,1}$ to $X$.}
of the supergroup $\mathbb{A}^{0,1}\rtimes(\mathbb{A}^1)^\times$
defined as follows. Let $o\in \lb^{1}=\mathbb{A}^{0,1}(\lb)$ we need
to define the corresponding automorphism of $\Pi
TX(\lb)=X(\lb[\xi])$.  We accomplish that by defining a morphism in
$\Lambda$ from $\lb[\xi]$ to itself via
$\xymatrix{\lb\ar[r]^{Id}&\lb}$ and $\xi\mapsto\xi+o$, this induces
the required automorphism. Similarly we define the automorphism
corresponding to $e\in(\lb^{0})^\times=(\mathbb{A}^1)^\times(\lb)$
by defining a morphism in $\Lambda$ via
$\xymatrix{\lb\ar[r]^{Id}&\lb}$ and $\xi\mapsto e\xi$.

At this point, for the sake of concreteness, let us assume that $X$
is prorepresentable (in the sense of Definition \ref{prorepdef}).
This is always the case in our setting.  We are now ready to define
the differential graded sheaf $\Omega_{X/\zp}$ of $\zp$-modules on
$[X]$ the body of $X$. Its hypercohomology will be called the de
Rham cohomology of $X$. We denote it by $DR_{\zp}(X)$.

Let $U\subset [X]$ be an open subvariety, consider the $\zp$-algebra
of functions on $\Pi T(X|_U)$ (i.e. natural transformations to
$\mathbb{A}^{1,1}$).  More concisely we have the following.

\begin{definition}
Define the pre-sheaf $\Omega_{X/\zp}$ on $[X]$ by setting
$$\Omega_{X/\zp}=[\pi]_*\mathcal{O}_{\Pi TX}.$$
\end{definition}

Note that $\Omega_{X/\zp}$ carries a grading induced by the action
of $(\mathbb{A}^1)^\times$.  The sections of $\Omega_{X/\zp}^n$
should be thought of as differential $n$-forms on $X$.  We also have
a differential\footnote{Speaking informally, it is induced by the
infinitesimal action of $\A^{0,1}$.} $d$ (an operator which comes
from the canonical odd element\footnote{In fact $\eta$ is the
canonical element in the odd Lie algebra $\A^{0,1}(\zp[\eta])$ of
$\A^{0,1}$.} $\eta\in\A^{0,1}(\lb[\eta])$ for every $\lb$). More
precisely, suppose that $\varphi$ is a function on $\Pi TX$. Recall
that this means that for every $\lb$ there is a natural map
$\varphi:\Pi TX(\lb)\rightarrow \lb$. Define the function $d\varphi$
by specifying that for every $\lb$, it is the following composition
of maps:$$\Pi TX(\lb)\stackrel{i}{\rightarrow}\Pi
TX(\lb[\eta])\stackrel{a_\eta}{\rightarrow}\Pi
TX(\lb[\eta])\stackrel{\varphi}{\rightarrow}\lb[\eta]\stackrel{c}{\rightarrow}\lb$$
where $i$ is induced by the inclusion $\lb\subset\lb[\eta]$,
$a_\eta$ is the endomorphism of $\Pi
TX(\lb[\eta])=X(\lb[\eta][\xi])$ induced by the endomorphism of
$\lb[\eta][\xi]$ given by $\lb\stackrel{Id}{\rightarrow}\lb$,
$\eta\mapsto\eta$, $\xi\mapsto\xi+\eta$; $\varphi$ is self
explanatory and $c$ reads off the coefficient of $\eta$.

One readily checks that $d$ increases the degree by
one.\footnote{The connection between these abstract definitions and
the usual de Rham complex is made explicit in Lemma \ref{localiso}.}

\begin{definition}\label{derhamcoh}
Let $X$ be a $p$-adic superspace. We define the de Rham cohomology
of $X$ as the hypercohomology of $\Omega^\bullet_{X/\zp}$, i.e.,
$$DR_{\zp}(X)=\mathbb{H}([X],\Omega^\bullet_{X/\zp}).$$
\end{definition}

Note that it is easy to see from the definitions that $DR_{\zp}(-)$
is a contravariant functor from the category of superspaces to the
category of graded $\zp$-modules. Thus any endomorphism of $X$
induces an endomorphism of $DR_{\zp}(X)$.

\noindent\textbf{Example.} Let us apply the definitions in the
simple example of a line. In this case $X(\lb)=\lb^0$ and $\Pi
TX(\lb)=X(\lb[\xi])=(\lb[\xi])^0=\lb^0+\lb^1\xi$.  Thus $\Pi
TX=\mathbb{A}^{1,1}$ as expected.  The body of $X$ is affine, so to
compute $DR_{\zp}(X)$ we need only compute the cohomology of the
complex of $\zp$-modules $\Gamma(\mathcal{O}_{\mathbb{A}^{1,1}})$.
By Theorem \ref{theorem:ankm} we know that as a $\zp$-module it is
$S\oplus S\xi$ where $S=\{\sum a_i x^i|a_i\in\zp, a_i\rightarrow
0\}$.  The reader is invited to verify that the action of
$(\mathbb{A}^1)^\times$ puts $S$ in degree $0$ and $S\xi$ in degree
$1$, while the differential acts by $\xi\partial_x$.  Thus
$DR_{\zp}(X)=H_{dR}(\widehat{\mathbb{A}^1}_p)$.

More generally, by unraveling the definitions we obtain the
following two Lemmas that bridge the gap between the $p$-adic
superspace approach and the usual situation.

\begin{lemma}\label{localiso}
Let $\text{Spec}A$ be a smooth variety over $\zp$, then
$$DR_{\zp}(X_{\text{Spec}A})=H_{dR}(\widehat{\text{Spec}A}_p).$$
\end{lemma}
\begin{proof}
Let $X=X_{\text{Spec}A}$, since $[X]$ is affine the left hand side
is computed by the complex of global functions on $\Pi TX$.

But $\Pi
TX(\lb)=X(\lb[\xi])=Hom_{Super}(A,\lb[\xi])=Hom_{Super}(\Omega_A,\lb)$.
Recall that $\Omega_A$ is the supercommutative ring generated by $a$
(even) and $da$ (odd) for $a\in A$ subject to the usual relations
(most important being the Leibniz Rule).

So $\Pi TX$ is prorepresented by $\Omega_A$, thus $\mathcal{O}_{\Pi
TX}=\widehat{(\Omega_A)}_p$ (just like in Lemma \ref{difference}).
This is exactly the complex that computes the right hand side, but
we still need to verify that this identification is compatible with
the differentials.  It is sufficient to check the compatibility with
the actions of $\mathbb{A}^{0,1}\rtimes(\mathbb{A}^1)^\times$.

Note that the action of $\A^{0,1}$ on $\Omega_A$ is given by the
coaction (which is an algebra morphism)
$$\Omega_A\rightarrow \Omega_A[\xi]$$ $$a\mapsto a+da\xi,\quad da\mapsto
da$$ and the action of $(\mathbb{A}^1)^\times$ by
$$\Omega_A\rightarrow \Omega_A[x^{\pm 1}]$$ $$a\mapsto a,\quad da\mapsto da
x$$ while $$Hom_{Super}(\Omega_A,\lb)=Hom_{Super}(A,\lb[\xi])$$
$$\{a\mapsto f(a), da\mapsto \varphi(a)\} \leftrightarrow\{a\mapsto
f(a)+\varphi(a)\xi\}.$$  It is now straightforward to check that the
action on $\Pi TX$ is obtained in this case from the one on
$\Omega_A$.

\end{proof}

\begin{lemma}\label{bridge}
Let $V$ be a smooth variety over $\zp$. Then
$$DR_{\zp}(X_V)=H_{dR}(\hat{V}_p).$$
\end{lemma}

\begin{proof}
The left hand side is by definition the hypercohomology of
$\Omega_{X_V/\zp}^\bullet$ on $[X_V]$, while the right hand side is
the hypercohomology of $\widehat{(\Omega_{V/\zp}^\bullet)}_p$ on
$V|_{\fp}$.  The two spaces $[X_V]$ and $V|_{\fp}$ are the same
(Lemma \ref{difference}), so the question is about comparing the two
complexes of sheaves locally.  They are the same by the proof of
Lemma \ref{localiso}.
\end{proof}

\subsection{Functions on $\widetilde{pt}$.}\label{pita}
In this section we study the most basic and at the same time the
most crucial example of a DP neighborhood, namely that of a point in
a line.

\begin{definition}
Denote by $\widetilde{pt}$ the infinitesimal neighborhood of the
origin in $\A^1$.
\end{definition}

In this section we are concerned with describing explicitly the
functions on $\widetilde{pt}$. This is the key step in the
subsequent cohomology computations.  A more category theory minded
reader may wish to visit the Appendix before going any further.

One sees immediately from the definitions that
$$\widetilde{pt}(\lb)=pB+\lb^{0+}.$$ This is our old friend $\A^{0|1,0}$,
and has a subfunctor that we will denote by $\widetilde{pt}_{0}$, it
is defined by $$\widetilde{pt}_{0}(\lb)=\lb^{0+}.$$  It is the
functions on $\widetilde{pt}_{0}$, i.e natural transformations to
$\A^1$ that we describe first.\footnote{Here we do not need
$\A^{1,1}$ as everything is purely even.} Let $f$ be one such
transformation, our intention is to show that for $w\in\lb^{0+}$, we
have that $f(w)=\sum_{i=0}^{\infty}a_i w^i/i!$ with $a_i\in \zp$
determining $f$.

\begin{remark}
It is clear that any expression $\sum_{i=0}^{\infty}a_i w^i/i!$
gives a function as it can be evaluated at any element of
$pB+\lb^{0+}$, i.e. we do have a map from such expressions to
functions. However the injectivity and surjectivity of this map
remains to be demonstrated below.
\end{remark}

First we need a Lemma.

\begin{lemma}\label{functionsonw}
Let
$w=\xi_{j_1}\xi_{j_2}+...+\xi_{j_{2k-1}}\xi_{j_{2k}}\in\Lambda_{\zp/p^N\zp}^{0+}$,
then $$f(w)=\sum_{i=0}^{k}a_i w^i/i!$$ and $a_i\in \zp/p^N\zp$
depend only on $f$.\footnote{The $a_i$ are defined inductively in
the proof.}
\end{lemma}
\begin{proof}
We proceed by induction on $k$.  If $k=0$ then $w=0$ and so by
functoriality of $f$, $f(w)\in \zp/p^N\zp$, define $a_0$ to be
$f(w)$ and we are done.

Assume that the Lemma is true for $k\leq n$.  Let $k=n+1$,
$$w=\xi_{j_1}\xi_{j_2}+...+\xi_{j_{2k-1}}\xi_{j_{2k}}=:x_1+...+x_k,$$
and setting $f(w)=\sum_I a_I\xi_I$ consider
$\xi_I=\xi_{j_1}...\xi_{j_{2i}}$.  Note that by functoriality
$i\leq k$, and if $i=k$ then there is only one such $\xi_I$,
denote its coefficient by $a_i$ (we have now defined $a_{n+1}$).
If $i<k$ define a map $\phi$ from $\Lambda_{\zp/p^N\zp}$ to itself
by sending $\xi_{j_{s}}$ to $\xi_s$ and the rest of $\xi$'s to
$0$.

If $\xi_I=x_{s_1}...x_{s_i}$ then
$\phi(w)=\xi_1\xi_2+...+\xi_{2i-1}\xi_{2i}$, so
$f(\phi(w))=...+a_i \xi_1\xi_2...\xi_{2i-1}\xi_{2i}$ by the
induction hypothesis and $\phi(f(w))=...+a_I
\xi_1\xi_2...\xi_{2i-1}\xi_{2i}$, so that $a_I=a_i$.

If $\xi_I\neq x_{s_1}...x_{s_i}$ then $\phi(w)$ has fewer than $i$
summands yet is of the form $\xi\xi+...+\xi\xi$ so that we may use
the induction hypothesis to conclude that the top degree of
$f(\phi(w))$ is less than $2i$ whereas
$\phi(\xi_I)=\xi_1\xi_2...\xi_{2i-1}\xi_{2i}$ has degree $2i$, so
that $a_I=0$.

So $f(w)=\sum_I a_I\xi_I=\sum_i a_i
x_{s_1}...x_{s_i}=\sum_{i=0}^{n+1}a_i w^i/i!$, and we are almost
done.  Namely, we demonstrated that any function $f$, when
restricted to $w\in\Lambda_{\zp/p^N\zp}^{0+}$ can be written as a DP
polynomial with coefficients in $\zp/p^N\zp$ of degree at most $k$.
However, it is immediate that such an expression is unique, since
$w^i/i!$ for $i=0,...,k$ form a basis of the free
$\zp/p^N\zp$-submodule of $\Lambda_{\zp/p^N\zp}^{0}$ that they span.

\end{proof}

By functoriality of $f$ we obtain, by considering the above Lemma
for all $N$, that the coefficients $a_i$ are given for all $N$ by
the images under the natural projection of $a_i\in\zp$.

\begin{theorem}\label{functions}
Let $e\in\lb^{0+}$, then $f(e)=\sum_{i=0}^{\infty}a_i e^i/i!$, with
$a_i$ as above.
\end{theorem}

\begin{proof}
Let $e\in\lb^{0+}$, $e=\sum_i b_i e_{i1} e_{i2}$ where $e_{ij}$ are
odd monomials in $\xi$'s.  Let $p^N B=0$.  Define a map $\varphi$
from $\Lambda_{\zp/p^N\zp}$ to $\lb$ by $\zp/p^N\zp\rightarrow B$
being the structure morphism, and $\xi_{2i-1}\mapsto b_i e_{i1}$ and
$\xi_{2i}\mapsto e_{i2}$, so that $w=\sum_i
\xi_{2i-1}\xi_{2i}\mapsto e$.  So
$f(e)=f(\varphi(w))=\varphi(f(w))=\varphi(\sum_{i=0}^{length(w)}a_i
w^i/i!)=\sum_{i=0}^{length(e)}a_i
\varphi(w)^i/i!=\sum_{i=0}^{length(e)}a_i
e^i/i!=\sum_{i=0}^{\infty}a_i e^i/i!$. Here the $length$ of an
element in $\lb^+$ is the minimal number of monomials (in the odd
variables) that are needed to write it down.
\end{proof}

\begin{remark}
A consequence of this result is that while our choice of $w$ in
Lemma \ref{functionsonw} is fairly arbitrary, for instance one can
reorder the coordinates, this does not in any way affect the
coefficients $a_i$.
\end{remark}

Now let us consider the functor $\widetilde{pt}(\lb)=pB+\lb^{0+}$
itself. We claim that the functions are still of the form
$\sum_{i=0}^{\infty}a_i x^i/i!$ with coefficients in
$\zp$.\footnote{Here in particular we must assume that $p>2$
otherwise this expression does not define a function in general
and the situation becomes more complicated.} We reduce to the
previous case to prove the following lemma, from which the claim
follows easily.

\begin{lemma}\label{lem}
Let $A=\zp/p^N\zp[x]$ and
$$w=px+\xi_{1}\xi_{2}+...+\xi_{{2k-1}}\xi_{{2k}}\in
pA+\Lambda_A^{0+},$$ then $$f(w)=\sum_{i=0}^{\infty}a_i w^i/i!$$
and $a_i\in\zp$ depend only on $f$.
\end{lemma}
\begin{proof}
First we need to define $a_i\in\zp$.  Recall the subfunctor
$\widetilde{pt}_{0}$ of $\widetilde{pt}$ that sends $\lb$ to
$\lb^{0+}$.  If we restrict $\widetilde{pt}_{0}$ to the subcategory
$\Lambda^N$ of Grassmann rings with coefficients in $B$, with $p^N
B=0$, then by Theorem \ref{functions}, $f|_{\widetilde{pt}_{0}}$
determines (and is determined on $\Lambda^N$ by)
$\{a_i^N\in\zp/p^N\zp\}$.  We observe that by functoriality of $f$
we may take the inverse limit over $N$ to obtain $\{a_i\in\zp\}$
that determine $f|_{\widetilde{pt}_{0}}$ on $\Lambda$.  Let
$\widetilde{f}$ be a new function on $\widetilde{pt}$ defined by
$\widetilde{f}(e)=\sum a_i e^i/i!$ for $e\in pB+\lb^{0+}$, so that
$\widetilde{f}$ agrees with $f$ on $\widetilde{pt}_{0}$.  We want to
show that they agree on $w$ also.

For any $N$ and $n$, let us define a map $\phi$ from $\Lambda_A$
to $\Lambda_{\zp/p^N\zp}$ by $\xi_i\mapsto\xi_i$ and
$x\mapsto\eta_1\eta_2+...+\eta_{2n-1}\eta_{2n}$.  Note that under
this map $$w\mapsto
p(\eta_1\eta_2+...+\eta_{2n-1}\eta_{2n})+\xi_{1}\xi_{2}+...+
\xi_{{2k-1}}\xi_{{2k}}\in\Lambda_{\zp/p^N\zp}^{0+}$$ and so
$$\phi(f(w))=f(\phi(w))=\sum a_i
(\phi(w))^i/i!=\phi(\widetilde{f}(w)).$$  Note that setting
$f(w)=\sum c_i^N x^i$ with $c_i^N\in \zp/p^N\zp[\xi_j]$,
functoriality implies that we have $c_i\in \zp[\xi_j]$ such that
$f(w)=\sum c_i x^i$ for all $N$.

To show that $f(w)=\widetilde{f}(w)$ it suffices to consider the
following situation.  Let $b_i\in\zp[\xi_j]$, define $g=\sum b_i
x^i\in\Lambda_A$, suppose that $$0=\phi(g)=\sum b_i
(\eta_1\eta_2+...+\eta_{2n-1}\eta_{2n})^i$$ for all $n$ and $N$.
Since if $i\leq n$ then $0=\phi(g)$ implies that $i!b_i=0$ we see
that $b_i=0$ in $\zp[\xi_j]$ and $g=0$.  Take
$g=f(w)-\widetilde{f}(w)$ and we are done.
\end{proof}

\begin{theorem}\label{functionsplus}
Let $e\in pB+\lb^{0+}$, then $f(e)=\sum_{i=0}^{\infty}a_i e^i/i!$.
\end{theorem}

\begin{proof}
Let $e\in pB+\lb^{0+}$, $e=pb+\sum_i b_i e_{i1} e_{i2}$ where
$e_{ij}$ are odd monomials in $\xi$'s and $b\in B$. Suppose that
$p^N B=0$. Define a map $\varphi$ from $\Lambda_A$ to $\lb$ by
$x\mapsto b$, $\xi_{2i-1}\mapsto b_i e_{i1}$ and $\xi_{2i}\mapsto
e_{i2}$, so that $w=px+\sum_i \xi_{2i-1}\xi_{2i}\mapsto e$.  So
that
$$f(e)=f(\varphi(w))=\varphi(f(w))=\varphi(\sum a_i w^i/i!)=\sum
a_i \varphi(w)^i/i!=\sum a_i e^i/i!.$$
\end{proof}

\begin{remark}
To summarize the above, we have a canonical map from
$\zp\left<x\right>$ to functions on $\widetilde{pt}$.  This map, as
is explicitly shown in the Lemmas above, is surjective.  The fact
that it is injective follows from the observation at the end of the
proof of Lemma \ref{functionsonw}.

\end{remark}

\subsection{The case of $p=2$.}\label{evenprime}

As mentioned before the case of the even prime does not fit into the
framework described.  The issue is that $\sum_{i=0}^{\infty}p^i/i!$
is convergent in $\zp$ only for $p>2$. It follows that for the case
$p=2$, the functions on the Grassmann neighborhood of a point in the
line are not simply DP power series with coefficients in $\zp$,
rather they form a subset of these with some conditions on the
coefficients.  While it is possible to describe them explicitly one
immediately sees that the homotopy of Lemma \ref{homotopy} no longer
exists.  Consequently one can not prove the cohomology invariance of
Grassmann thickening.

It seems one can introduce an alternate framework that works for all
primes $p$.  We briefly outline it here. The idea is to introduce
$\hat{\Lambda}$, an enlargement of our main category $\Lambda$ which
includes Grassmann rings with an infinite number of variables that
allow certain infinite sums as elements. More precisely, we consider
rings $\lb=B[\xi_1,\xi_2,...]$ where elements have the form $\sum
b_i w_i$ where $b_i\in B$ and $w_i$ are monomials in $\xi$'s of
degree at most $N$ where $N$ is fixed for each element.\footnote{We
still require that $p^M B=0$ for $M$ large, thus the canonical DP
ideal is still nilpotent, but it does not have to be DP nilpotent.} Thus
$\sum \xi_{2i-1}\xi_{2i}$ is an element, while
$\sum_{i=1}^{\infty}\prod_{j=1}^i \xi_j$ is not.

One does not get the same functions as before for the case of the
Grassmann neighborhood of a point in the line\footnote{Instead of
power series with DP one gets an extra condition that the
coefficients tend to $0$ in $\zp$.}, but the homotopy of Lemma
\ref{homotopy} now makes sense and so we can again show the
cohomology invariance of Grassmann thickening by modifying all of
the arguments accordingly (some of them simplify somewhat).

Finally, note that the very definition of the Grassmann algebra
needs modification by the addition of an extra axiom that
$\xi^2=0$ for $\xi$ odd.

\subsection{De Rham cohomology in the smooth case, a comparison.}\label{smooth}
Recall that to a variety $V$ over $\zp$ , considered as  a functor from
$\zp$-algebras to sets, we can associate a superspace $X_V$ with
$X_V(\lb)=V(\lb^{0})$.  If $V$ is smooth, then we may consider the
usual de Rham cohomology of $V$ and compare it to the
$DR_{\zp}(X_V)$.  In general the two are not the same, however if
$V$ is projective then they are isomorphic.

\begin{lemma}\label{derhamtodr}
Let $V$ be a smooth projective variety over $\zp$ then
$$H_{dR}(V)\simeq DR_{\zp}(X_V).$$
\end{lemma}

\begin{proof}
By Lemma \ref{bridge} we need only compare $H_{dR}(\hat{V}_p)$ with
$H_{dR}(V)$.  The fact that they are isomorphic in the projective
case was pointed out to us by A. Ogus, and we provide a sketch of a
proof. For the relevant facts about formal schemes and the theorem
on formal functions we refer to \cite{hartshorne}.

By definition, the de Rham cohomology $H_{dR}(V)$ is computed as the
hypercohomology of the complex $\Omega_{V/\zp}^\bullet$, which can
be obtained as the cohomology of the total complex associated to the
double complex of $\zp$-modules $\oplus_I
\Omega_{V/\zp}^\bullet(U_I)$, where $I$s are finite subsets
$\{i_1,...,i_s\}$ of $\{0,...,n\}$, $U_I=U_{i_1}\cap...\cap U_{i_s}$
and $U_0,...,U_n$ form an open affine cover of $V$.  The horizontal
differentials are de Rham differentials and the vertical ones are
\v{C}ech differentials.

Clearly we have a morphism of double complexes $$\alpha:\oplus_I
\Omega_{V/\zp}^\bullet(U_I)\rightarrow \oplus_I
\widehat{(\Omega_{V/\zp}^\bullet(U_I))}_p$$ and the double complex
on the right computes $H_{dR}(\hat{V}_p)$.

The above morphism on the $E_1$ term becomes
$$\alpha:H^\bullet(\Omega_{V/\zp}^\bullet)\rightarrow
H^\bullet(\widehat{(\Omega_{V/\zp}^\bullet)}_p)$$ and it factors as
follows $$H^\bullet(\Omega_{V/\zp}^\bullet)\rightarrow
\widehat{H^\bullet(\Omega_{V/\zp}^\bullet)}_p\rightarrow
H^\bullet(\widehat{(\Omega_{V/\zp}^\bullet)}_p).$$ Recall that
projective morphisms preserve coherence and so
$H^i(\Omega_{V/\zp}^j)$ is a finitely generated $\zp$-module.
Because $\zp$ is complete, by Theorem 9.7 in \cite{hartshorne}, for
example, we have that the first arrow is an isomorphism.  The second
arrow is an isomorphism by the theorem on formal functions.

Since $\alpha$ is an isomorphism on $E_1$, it induces an isomorphism
$$\alpha:H_{dR}(V)\stackrel{\sim}{\rightarrow}H_{dR}(\hat{V}_p).$$

\end{proof}

Given a smooth projective $V$ over $\zp$ we would like to define the
action of $Fr$ on its de Rham cohomology.  By above it suffices to
do so for $DR_{\zp}(X_V)$.  As explained in Sec.
\ref{supspnhoodsandfrob}, we have an action of $Fr$ on
$\widetilde{X_V}$ (the neighborhood of $X_V$ in $\mathbb{P}^n$) and
so on $DR_{\zp}(\widetilde{X_V})$.  Showing that
$DR_{\zp}(\widetilde{X_V})$ is isomorphic to $DR_{\zp}(X_V)$ would
accomplish our goal.

\begin{remark}

Another consequence of this isomorphism is that the de Rham
cohomology of $V$ depends only on $V|_{\fp}$ because that is true of
$\widetilde{X_V}$.  This means in particular that for $W$ smooth
projective over $\fp$, the de Rham cohomology of $\widetilde{W}$
(its Grassmann neighborhood in the projective space over $\zp$)
coincides with the usual crystalline cohomology of $W$.  That is we
give a super-geometric construction of the DP envelope of $W$ in
$\mathbb{P}^n_{\zp}$.

\end{remark}

Observe that we have $i:X_V\hookrightarrow\widetilde{X_V}$ thus
also a natural map $$i^*:DR_{\zp}(\widetilde{X_V})\rightarrow
DR_{\zp}(X_V)$$ that we will show is an isomorphism.  It suffices
to prove that $$i^*:\Omega_{\widetilde{X_V}/\zp}\rightarrow
\Omega_{X_V/\zp}$$ is a quasi-isomorphism of sheaves on
$V|_{\fp}$. Thus the question becomes local and we may assume,
after induction on the codimension, that the situation is as
follows.

Let $U\subset U'$ be a pair of smooth affine varieties such that $U$
is cut out of $U'$ by a function $g$ on $U'$. In this case we will show that
\begin{align*}
\widetilde{X}_U=X_U\times\widetilde{pt}
\end{align*}
i.e.  the infinitesimal neighborhood of $U$ in $U'$ is a direct
product of the $p$-adic superspaces $X_U$ (associated to $U$) and
our $\widetilde{pt}$. Compare this with the Example following
Definition \ref{nhood} where this is discussed in the case when the
function $g$ is linear and $U'$ is an affine space. The general case
is demonstrated below. We will assume that $U'=\text{Spec} A$ and
$U=\text{Spec} A/g$, where $g\in A$. Then unraveling the definitions
we see that
\begin{align*}
\widetilde{X}_U(\lb)&=\{\Hom(A,\lb)|g\mapsto
pB+\lb^+\}\\
&=\bigcup_n\Hom(A/g^n,\lb)\\
&=\Hom^{cont}(\hat{A}_g,\lb)\\
&=\Hom^{cont}(A/g[[x]],\lb)\\
&=\Hom(A/g,\lb)\times(pB+\lb^{0+})\\
&=X_U(\lb)\times\widetilde{pt}(\lb)\\
\end{align*}  where $\Hom^{cont}$ denotes continuous homomorphisms ($\hat{A}_g$ is a topological ring and
$\lb$ is equipped with the discrete topology).  The inclusion
$X_U\subset \widetilde{X}_U$ is simply $X_U=X_U\times pt\subset
X_U\times\widetilde{pt}=\widetilde{X}_U$.

Next we show that the functions on $X_U\times\widetilde{pt}$ are
what was expected, namely if $R=\widehat{(A/g)}_p$ then:

\begin{theorem}\label{relfunctions} Any natural transformation $f$ from
$X_U\times\widetilde{pt}$ to $\mathbb{A}^1$ is given by
$$f(\phi,e)=\sum_{i=0}^{\infty}\phi(r_i)e^i/i!$$ for any $\phi\in
\Hom(A/g,\lb)$ and $e\in pB+\lb^{0+}$, where $r_i\in R$ depend
only on $f$.
\end{theorem}
\begin{remark}
The proof below remains valid for the case when $X$ is given by
$X(\lb)=\Hom(\zp\left<x_i\right>,\lb)$ i.e.
$R=\zp\left<x_i\right>$. This justifies the induction on the
codimension.
\end{remark}

\begin{proof}
Begin by noting that in the proofs of Lemma \ref{lem} and Theorem
\ref{functionsplus} we can replace $\zp$ by any $p$-adically
complete ring without zero divisors.  In particular these results
remain valid when $\zp$ is replaced by our $R=\widehat{(A/g)}_p$.

Thus let us define a new category $\Lambda(R)$ consisting of
Grassmann rings with coefficients in $R$-algebras with nilpotent
$p$-action. Denote by $\widetilde{pt}_R$ and $\mathbb{A}^1_R$ the
restrictions of similarly named functors from $\Lambda$ to
$\Lambda(R)$, so that they are now functors from $\Lambda(R)$ to
$Sets$. As before, if $f$ is a natural transformation from
$\widetilde{pt}_R$ to $\mathbb{A}^1_R$, then
$$f(e)=\sum_{i=0}^{\infty}r_i e^i/i!,$$ where $e\in pB+\lb^{0+}$,
and $r_i\in R$ depend only on $f$.

Observe that any natural transformation $f$ from
$X_U\times\widetilde{pt}$ to $\mathbb{A}^1$ defines
$\widetilde{f}:\widetilde{pt}_R\rightarrow\mathbb{A}^1_R$ by
$\widetilde{f}(e)=f(\phi,e)$, where $e\in
pB+\lb^{0+}\in\Lambda(R)$ and $\phi$ is the structure morphism. Of
course the $R$-module structure on $B$ provides us with a
$\zp$-morphism $\phi:R\rightarrow B\rightarrow \lb$, however as
$\phi$ factors through $R/p^N R=(A/g)/p^N (A/g)$ it determines a
unique morphism from $A/g$ to $\lb$ and so an element of
$X_U(\lb)$.  Conversely, any element of $X_U(\lb)$ that factors
though $B$ makes $\lb$ into an element of $\Lambda(R)$. By above
$$\widetilde{f}(e)=\sum_{i=0}^{\infty}r_i e^i/i!$$ that is
$$f(\phi,e)=\sum_{i=0}^{\infty}\phi(r_i)e^i/i!$$ for all
$\phi:A/g\rightarrow\lb$ that factor through $B$.

Let $(\phi,e)$ be an arbitrary element of
$X_U\times\widetilde{pt}(\lb)$, assume that $p^N B=0$ so that
$\phi$ factors through $R/p^N R$. Let $e=pb+\sum b_i e_{i1}
e_{i2}$ and proceed as in the proof of Theorem
\ref{functionsplus}. Define a morphism $\varphi$ from
$\Lambda_{(R/p^N R)[x]}$ to $\lb$ by
\begin{align*}
\phi:R/p^N R&\rightarrow\lb\\
x&\mapsto b\\
\xi_{2i-1}&\mapsto b_i e_{i1}\\
\xi_{2i}&\mapsto e_{i2}.
\end{align*}  Consider the
element $(\pi,w)\in X_U\times\widetilde{pt}(\Lambda_{(R/p^N
R)[x]})$ where $$\pi:A/g\rightarrow \Lambda_{(R/p^N R)[x]}$$ is
the projection onto $R/p^N R\subset \Lambda_{(R/p^N R)[x]}$ and
$$w=px+\sum_i \xi_{2i-1}\xi_{2i},$$ then
$\varphi(\pi,w)=(\phi,e)$.\footnote{One should really write
$X_U\times\widetilde{pt}(\varphi)((\pi,w))$, but that is too
cumbersome.} Thus
$$f(\phi,e)=f\varphi(\pi,w)=\varphi
f(\pi,w)=\varphi(\sum\pi(r_i)w^i/i!)=\sum\phi(r_i)e^i/i!.$$

\end{proof}

We are finally able to complete the proof of Theorem
\ref{theorem:ankm} which we restate as a Corollary below.

\begin{cor}\label{ankm}
The functions on $\mathbb{A}^{n|k,m}$ are isomorphic as a ring to
$$\left\{\sum_{I,J,T\geq 0}a_{I,J,T}x^I \xi^J y^T/T!\right\}$$ where $x_i$ and
$y_i$ are even and $\xi_i$ are odd, and $a_{I,J,T}\in\zp$ with
$a_{I,J,T}\rightarrow 0$ as $I\rightarrow \infty$.
\end{cor}

\begin{proof}
Using Theorems \ref{functionsplus} and \ref{relfunctions} with
induction we see that the ring of functions on
$\mathbb{A}^{0|k,0}$ is $\zp\left<y_1,...,y_k\right>$.  Gluing
this fact with the proven part of Theorem \ref{theorem:ankm} using
Theorem \ref{relfunctions} we obtain the desired result.
\end{proof}

Denoting the functions described in the Theorem \ref{relfunctions}
by $R\left<x\right>$ and observing that the functions on $X_U$ are
given by $R$ we are done by the following Lemma.

\begin{lemma}\label{homotopy}
The natural map $\pi:\Omega_{R\left<x\right>}\rightarrow\Omega_R$
is a quasi-isomorphism.
\end{lemma}

\begin{proof}
In fact we show that the equally natural map
$\rho:\Omega_R\rightarrow\Omega_{R\left<x\right>}$ is a homotopy
inverse.  Note that $\pi\circ \rho=\text{Id}_{\Omega_R}$, let
$F=\rho\circ \pi$, we must show that there is a homotopy $h$ such
that $\text{Id}_{\Omega_{R\left<x\right>}}-F=d\circ h + h\circ d$.

It follows immediately from the abstract definitions and by using
Theorem \ref{relfunctions}, that any
$w\in\Omega^{s}_{R\left<x\right>}$ can be written uniquely as
$$w=\sum_{i=0}^{\infty}\alpha_i x^i/i!+\sum_{i=0}^{\infty}\beta_i
x^i/i! dx$$ where $\alpha_i\in \Omega^s_R$ and
$\beta_i\in\Omega^{s-1}_R$.  Let
$$h(w)=(-1)^{s-1}\sum_{i=0}^{\infty}\beta_i x^{i+1}/(i+1)!,$$ then a
straightforward calculation shows that $h$ is the desired
homotopy.
\end{proof}

At this point we have proven:

\begin{lemma}\label{new}
Let $U\subset U'$ be a closed embedding of smooth varieties over $\zp$, then the natural map $$DR_{\zp}(\widetilde{X_U})\rightarrow
DR_{\zp}(X_U)$$ is an isomorphism.
\end{lemma}

And the main result of the section, namely:

\begin{theorem}\label{main}
Let $V$ be a smooth projective variety over $\zp$, let $X_V$ be the
associated $p$-adic superspace, and $\widetilde{X_V}$ the DP
neighborhood of $V$ inside $\mathbb{P}^n$, then $$H_{dR}(V)\simeq
DR_{\zp}(\widetilde{X_V}).$$
\end{theorem}

\begin{cor}\label{frobactoncoh}
Let $V$ be a smooth projective variety over $\zp$, then one has an
action of Frobenius on $H_{dR}(V)$.
\end{cor}
\begin{proof}
By Lemma \ref{frobenius} we have an action of Frobenius on
$\widetilde{X_V}$ and so on its de Rham cohomology, which is
isomorphic to $H_{dR}(V)$.
\end{proof}

\begin{theorem}\label{crysderham}
Let $W$ be a smooth projective variety over $\fp$, then the
crystalline cohomology of $W$ is isomorphic to the de Rham
cohomology of the $p$-adic superspace $\widetilde{W}$ (the Grassmann
neighborhood of $W$ in $\mathbb{P}^n_{\zp}$), i.e.,
$$H_{crys}(W)\simeq DR_{\zp}(\widetilde{W}).$$
\end{theorem}
\begin{proof}
Recall one of the key points of the theory of crystalline cohomology, namely, if $D$ is the $DP$ neighborhood of $W$ in $\mathbb{P}^n_{\zp}$, then $H_{crys}(W)$ is canonically isomorphic to $\mathbb{H}(W,\Omega^\bullet_{D/\zp})$.  The latter is almost the same as our definition of the de Rham cohomology of $\widetilde{W}$. In fact we have a map of complexes $\Omega^\bullet_{D/\zp}\rightarrow\Omega^\bullet_{\widetilde{W}/\zp}$ that we will show induces an isomorphism on hypercohomology.

In order to compute the hypercohomology, consider the double complex associated to the restriction to $W$ of the open cover $U_i$ of $\mathbb{P}^n_{\zp}$ by $\mathbb{A}^n_{\zp}$s, just as in the proof of Lemma \ref{derhamtodr}.  The map of complexes above induces a map of double complexes $$\oplus_I\Omega^\bullet_{D/\zp}(U_I)\rightarrow \oplus_I\Omega^\bullet_{\widetilde{W}/\zp}(U_I)$$ which reduces the problem to a local one by considering the $E_1$ term of the associated spectral sequence.  (Note that here the de Rham differential is applied first, while in the Lemma \ref{derhamtodr} the \v{C}ech differential was applied first.)

Thus it suffices to show that the map $H_{dR}(D_I)\rightarrow DR_{\zp}(\widetilde{W}_I)$ is an isomorphism.  Let $V_I$ be any smooth lifting of $W_I$ to $\zp$.  We then have a sequence of maps: $$H_{dR}(D_I)\rightarrow DR_{\zp}(\widetilde{W}_I)\rightarrow DR_{\zp}(X_{V_I}),$$ where the second map is induced by $X_{V_I}\rightarrow\widetilde{W}_I$ and is an isomorphism by Lemma \ref{new} since $\widetilde{W}_I=\widetilde{X_{V_I}}$.  Observe that the composition is the map $H_{dR}(D_I)\rightarrow H_{dR}(\widehat{(V_I)}_p)$, which is also an isomorphism by the $DP$ Poincar\'{e} Lemma.  And so we conclude that the first map is an isomorphism.

\end{proof}

It is worth noting that the homotopy inverse $\rho$ that was used in
the proof of Lemma \ref{homotopy} can not be realized (in general)
as a restriction of a global map $r^*:DR_{\zp}(X_V)\rightarrow
DR_{\zp}(\widetilde{X_V})$.\footnote{In contrast with $\pi$ which is
a restriction of a global map
$i^*:DR_{\zp}(\widetilde{X_V})\rightarrow DR_{\zp}(X_V)$.}
Geometrically speaking we may not in general have a global
retraction $r$ of $\widetilde{X_V}$ onto $X_V$, its existence would
ensure that $i^*:DR_{\zp}(\widetilde{X_V})\rightarrow DR_{\zp}(X_V)$
is an isomorphism of filtered modules with respect to the Hodge
filtration.  Consequently, the canonical lift of the Frobenius
morphism $Fr$ to $DR_{\zp}(X_V)$ would preserve the Hodge filtration
$F^{\bullet} DR_{\zp}(X_V)$.  Furthermore, consider the following
local computation.  Let $x$ be a local function on
$\widetilde{X_V}$, then $Fr(x)=x^p+py$, where $y$ is some other
local function, so that $$Fr:fdx_1...dx_s\mapsto p^s
Fr(f)(x_1^{p-1}dx_1+dy_1)...(x_s^{p-1}dx_s+dy_s),$$ i.e. under the
assumption that a global retraction exists
$$Fr(F^s DR_{\zp}(X_V))\subset p^s F^s DR_{\zp}(X_V).$$  Neither
the invariance of the Hodge filtration under $Fr$ nor the
$p$-divisibility estimate need hold in the absence of the global
retraction, in Sec. \ref{hodge} we obtain some weaker
$p$-divisibility estimates that hold in general.

\section{The Frobenius map and the Hodge
filtration.}\label{hodge}

In this section we essentially follow B. Mazur\cite{mazur} with some
differences in the point of view (that is we find it more conceptual
to think of DP ideals and their DP powers).  We begin by recalling a
definition.

\begin{definition}
Let $I$ in $A$ be a DP ideal, then the $n$th DP power of $I$,
denoted $I^{[n]}$, is the ideal generated by the products
$x_1^{n_1}/n_1 !...x_k^{n_k}/n_k !$ with $x_i\in I$ and $\sum
n_i\geq n$.
\end{definition}

We point out that the Hodge filtration on the de Rham cohomology of
$X$ is simply the filtration on the functions on $\Pi TX$ given by
the DP-powers of the DP ideal $I_X$ of $X$ in $\Pi TX$.  More
precisely:

\begin{definition}
For a $p$-adic superspace $X$, define a filtration, $F_H^\bullet$ on
$\Omega_{X/\zp}$ by setting $F_H^i \Omega_{X/\zp}=I_X^{[i]}$.  This
filtration descends to $DR_{\zp}(X)$\footnote{The action of
$\A^{0,1}$ on $\Pi TX$ preserves $X$, thus the differential
preserves $I_X$ and its DP powers.} and let us still denote it by
$F_H^\bullet$ there.
\end{definition}

The following Lemma is immediate.

\begin{lemma}
Let $V$ be a smooth projective variety over $\zp$, then the
isomorphism $$H_{dR}(V)\simeq DR_{\zp}(X_V)$$ is compatible with the
Hodge filtration on the left and $F_H^\bullet$ on the right.
\end{lemma}

\begin{remark}
Because of the above Lemma we will use the notation $F_H^\bullet$ to
denote also the Hodge filtration on $H_{dR}(V)$.
\end{remark}

However the Frobenius map that we are interested in is defined on
$DR_{\zp}(\widetilde{X_V})$ not $DR_{\zp}(X_V)$. And while the two
are isomorphic as shown previously, this isomorphism is not
compatible with $F_H^\bullet$.  To fix this problem we proceed as
follows: replace the $F_H^\bullet$ filtration on
$DR_{\zp}(\widetilde{X_V})$ which is given by the ideal of
$\widetilde{X}$ in $\Pi T\widetilde{X}$ with the one given by the
DP-powers of the ideal of $X$ itself in $\Pi T\widetilde{X}$.

\begin{definition}
Let $X\subset\widetilde{X}\subset\Pi T\widetilde{X}$ be as above,
define a filtration $F_{DP}^\bullet$ on $DR_{\zp}(\widetilde{X})$ as
the induced filtration from $\Omega_{\widetilde{X}/\zp}$, where
$$F_{DP}^i \Omega_{\widetilde{X}/\zp}=I_X^{[i]}.$$
\end{definition}

In the particular case, namely the setting of Lemma \ref{homotopy}
(that is the key step in proving the general case), the definition
above becomes: $F_{DP}^\bullet\Omega_{R\left<x\right>}$ is defined
by $w\in F^s \Omega_{R\left<x\right>}$ if
$$w=\sum_{i=0}^{\infty}\alpha_i x^i/i!+\sum_{i=0}^{\infty}\beta_i
x^i/i! dx$$ where $\alpha_i\in \Omega^{\geqslant s-i}_R$ and
$\beta_i\in\Omega^{\geqslant s-1-i}_R$. It is then not hard to show
(using the observation that the homotopy of Lemma \ref{homotopy}
preserves the new filtration) that
$i^*:DR_{\zp}(\widetilde{X_V})\rightarrow DR_{\zp}(X_V)$ is an
isomorphism of filtered modules where $DR_{\zp}(\widetilde{X_V})$ is
endowed with the new filtration $F_{DP}^\bullet$ and $DR_{\zp}(X_V)$
has the old filtration $F_H^\bullet$.  Thus we have a refinement of
Theorem \ref{main}.

\begin{theorem}\label{mainprime}
Let $V$ be a smooth projective variety over $\zp$, let $X_V$ be the
associated $p$-adic superspace, and $\widetilde{X_V}$ the DP
neighborhood of $V$ inside $\mathbb{P}^n$, then $$H_{dR}(V)\simeq
DR_{\zp}(\widetilde{X_V})$$ is an isomorphism of filtered modules
with the Hodge filtration on the left and $F_{DP}^\bullet$ on the
right.
\end{theorem}

Though $Fr$ does not preserve it, $F_{DP}^\bullet$ is useful in the
computation of divisibility estimates.  Notice that
$Fr(I^{[k]}_X)\subset p^{[k]} \mathcal{O}_{\widetilde{X}}$
\footnote{A very useful formula to keep in mind is $\text{ord}_p
n!=\displaystyle\sum_{i=1}^{\infty}\left[\dfrac{n}{p^i}\right]
=\dfrac{n-S(n)}{p-1}$, where $S(n)$ is the sum of digits in the
p-adic expansion of $n$.}, where $I^{[k]}_X$ is the $k$-th DP power
of the ideal (also denoted by $I_X$) of $X$ in $\widetilde{X}$ and
$[k]=\text{min}_{n\geqslant k}\text{ord}_p (\frac{p^n}{n!})$ (thus
$(p)^{[k]}=(p^{[k]})$). Recalling the discussion at the end of Sec.
\ref{smooth}, we see that $$Fr(F_{DP}^s
DR_{\zp}(\widetilde{X_V}))\subset
p^{[s]}DR_{\zp}(\widetilde{X_V}).$$  In particular if $p>\text{dim}
(V)$ then the square brackets can be removed from $s$ and we obtain
the following Lemma.

\begin{lemma}\label{divest}
Let $V$ be a smooth projective variety over $\zp$ such that
$p>\text{dim} (V)$, then $$Fr(F^s_H H_{dR}(V))\subset p^s
H_{dR}(V).$$
\end{lemma}

A slightly finer statement can be derived from the above
observations, one actually has
\begin{align*}
Fr(F_{DP}^s DR_{\zp}(\widetilde{X_V}))&\subset \sum_{j<s} p^{[s-j]}
Fr(F_{DP}^{j}DR_{\zp}(\widetilde{X_V})) +p^s DR_{\zp}(\widetilde{X_V})\\
&\subset p Fr (DR_{\zp}(\widetilde{X_V}))+p^s DR_{\zp}(\widetilde{X_V}).\\
\end{align*} The latter was sufficient for Mazur to establish a conjecture of
Katz.

By analogous reasoning one can introduce new filtrations on the
cohomologies of $X$ and $\widetilde{X}$ by considering the DP ideal
of $X|_{\fp}$ in $\Pi TX$ and the DP ideal of $X|_{\fp}$ in $\Pi
T\widetilde{X}$.  The canonical isomorphism is now an isomorphism of
filtered modules with respect to these new filtrations and they are
preserved by $Fr$.  The new filtration on $DR_{\zp}(X_V)$ contains
the Hodge filtration and satisfies the same divisibility conditions.
In fact it can be easily described as follows (let us work with
$H_{dR}(V)$ since it is the same as $DR_{\zp}(X_V)$).  Let
$F_H^\bullet$ denote the usual Hodge filtration on $H_{dR}(V)$, then
the new filtration $F_N^\bullet$ can be described thus:
$F_N^{n}H_{dR}(V)=\sum_{s+t\geqslant n} p^{[s]}F_H^{t}H_{dR}(V)$.

\section{Appendix.}\label{colimit}
We investigate a property of a $p$-adic superspace that we call
prorepresentability.  It allows us to pass from particular examples
that we considered in this paper (namely $p$-adic superspaces that
arise in dealing with usual varieties over $\zp$) to a more general
class of $p$-adic superspaces that nevertheless share a lot of
properties with our examples.

Recall that we have defined a map of $p$-adic superspaces as a map
of the defining functors.  It is well known (Yoneda Lemma) that
the set $\text{Hom}(F,G)$ of natural transformations from a
functor $F$ to a functor $G$ can be easily calculated\footnote{The
source category does not matter and the target category is
$Sets$.} if $F$ is representable, i.e. $F$ is isomorphic to the
functor $h_A$, where
$$h_A(X)=\text{Hom}(A,X).$$  Namely, in this case we have
$$\text{Hom}(F,G)=G(A).$$

\begin{definition}
We say that a functor $F$ is prorepresentable if it is isomorphic
to a colimit of representable functors:$$F=\varinjlim_{A\in D}
h_A$$ for some diagram $D$.
\end{definition}

Thus $$\text{Hom}(F,G)=\varprojlim_{A\in D} G(A),$$ i.e. is a limit
of the sets $G(A)$, by Yoneda Lemma and continuity of
$\text{Hom}(-,-)$.  Notice that in the above definitions we use a
general definition of limits and colimits, a concise reference is
\cite{may}.  It is important to emphasize that the definition of a
prorepresentable functor in \cite{artin} is much more restrictive.

\begin{definition}\label{prorepdef}
We say that a $p$-adic superspace $X$ is prorepresentable, if its
defining functor is locally prorepresentable, i.e. $[X]$ has a cover
by open $U_i$ such that the defining functors of $X|_{U_i}$ are
prorepresentable.
\end{definition}

One can show that the property of a $p$-adic superspace $X$ being
prorepresentable is preserved by passing to the odd tangent space
$\Pi TX$.

\begin{theorem}
If a $p$-adic superspace $X$ is prorepresentable then so is $\Pi
TX$.
\end{theorem}
\begin{proof}
This statement is local, so we may assume that $X$ is
prorepresentable as a functor. Denote by $F_\xi$ the endo-functor of
$\Lambda$ that takes an object $A\in\Lambda$ to $A[\xi]$, i.e.
adjoins an odd variable. Note that $\Pi TX=X\circ F_\xi$.  Observe
that $F_\xi$ extends to $Super$. $Super$ is closed under limits and
$F_\xi$ commutes with limits since it has a left adjoint
$\Omega^\bullet_{-}$.  Thus we may assume that $X$ is representable
and since $\Omega^\bullet_{-}$ descends to $\Lambda$ we are done.
\end{proof}

One can show that all the $p$-adic superspaces we consider in this
paper are prorepresentable in our sense.  Here we give a detailed
proof of this fact for the most important functor
$\widetilde{pt}$, which is indeed prorepresentable (not just
locally prorepresentable). In a similar way one can prove the
local prorepresentability of functors corresponding to other
superspaces considered in the present paper. As an application we
show how to describe the functions on $\widetilde{pt}$ (Theorem
\ref{functionsplus}) using the above ideas.

Consider the ring $C[\xi_i]_{i=1}^{2n}$, i.e. a supercommutative
ring obtained from a commutative ring $C$ by adjoining $2n$ odd
anticommuting variables $\xi_j$. Let $G$ be the group acting on
$C[\xi_i]_{i=1}^{2n}$ generated by
$$\xi_{2k-1}\mapsto\xi_{2k},\,\xi_{2k}\mapsto-\xi_{2k-1}$$ and
$$\xi_{2k-1},\xi_{2k}\mapsto\xi_{2k'-1},\xi_{2k'}$$
$$\xi_{2k'-1},\xi_{2k'}\mapsto\xi_{2k-1},\xi_{2k}.$$

\begin{lemma}\label{ginvar}
The subring of $G$-invariants in $C[\xi_i]_{i=1}^{2n}$ is spanned by
$w_n^k/k!$ for $k=0,...,n$ where
$w_n=\xi_1\xi_2+...+\xi_{2n-1}\xi_{2n}$.  That is
$$(C[\xi_i]_{i=1}^{2n})^G=C\left<y\right>/(y)^{n+1}.$$
\end{lemma}
\begin{proof}
Proceed by induction on $n$.  For $n=0$ there is nothing to prove.
Let it be true for $n$.  Note that
\begin{align*}
C&[\xi_i]_{i=1}^{2(n+1)}\\
&=C[\xi_i]_{i=1}^{2n}\oplus C[\xi_i]_{i=1}^{2n}\xi_{2n+1} \oplus
C[\xi_i]_{i=1}^{2n}\xi_{2n+2}\oplus
C[\xi_i]_{i=1}^{2n}\xi_{2n+1}\xi_{2n+2}.\\
\end{align*}
Let $x\in C[\xi_i]_{i=1}^{2(n+1)}$ be $G$-invariant, then the
element of $G$ that ``switches" $\xi_{2n+1}$ and $\xi_{2n+2}$
ensures that
$$x\in C[\xi_i]_{i=1}^{2n}\oplus
C[\xi_i]_{i=1}^{2n}\xi_{2n+1}\xi_{2n+2}.$$

Considering the part of $G$ that acts on $C[\xi_i]_{i=1}^{2n}$ only
and using the induction hypothesis we see that $$x=\sum a_k w_n^k/k!
+ b_{k-1}w_n^{k-1}/(k-1)!\xi_{2n+1}\xi_{2n+2}.$$

Since the action of $G$ is degree preserving each homogeneous
component of $x$ is also $G$ invariant, thus $$a_k w_n^k/k! +
b_{k-1}w_n^{k-1}/(k-1)!\xi_{2n+1}\xi_{2n+2}$$ is $G$-invariant.

The element of $G$ that switches $\xi_1,\xi_2$ and
$\xi_{2n+1},\xi_{2n+2}$ ensures that $$a_k=b_{k-1},$$ so that
$$x=\sum a_k w_{n+1}^k/k!.$$
\end{proof}

Observe that if a group $G$ acts on a set $X$, then the fixed points
subset $X^G\subset X$ can be represented as a limit of the diagram
in $Sets$ consisting of two copies of $X$ and the arrows given by
the elements of $G$.  This is where Lemma \ref{ginvar} is used in
Theorem \ref{prorep} below.

\begin{theorem}\label{prorep}
The functor defining the $p$-adic superspace $\widetilde{pt}$ is
prorepresentable.
\end{theorem}
\begin{proof}
Consider the diagram in $\Lambda$ consisting of objects
$$\{(\zp/p^n\zp)[\xi_i]_{i=1}^{2m}|m,n\geq 0\}
\coprod \{(\zp/p^n\zp)[\xi_i]_{i=1}^{2m}|m,n\geq 0\}$$ with
morphisms between the copies given by elements of $G$, and the rest
of the morphisms given by the usual projections
$$(\zp/p^n\zp)[\xi_i]_{i=1}^{2m}\rightarrow(\zp/p^{n'}\zp)[\xi_i]_{i=1}^{2m}$$
for $n>n'$ and $$(\zp/p^n\zp)[\xi_i]_{i=1}^{2m}\rightarrow
(\zp/p^n\zp)[\xi_i]_{i=1}^{2m'}$$ for $m>m'$ mapping the extra
$\{\xi_j\}_{j>2m'}$ to $0$.

It follows from Lemma \ref{ginvar} that while the limit of the
above diagram does not exist in $\Lambda$, it exists in the
category $Super$, where $\Lambda$ is a full subcategory, and it is
equal to the ring with divided powers $\zp\left<x\right>$.
Correspondingly, after passing to the category of functors from
$\Lambda$ to $Sets$, the functor
$$\widetilde{pt}=\text{Hom}_{Super}(\zp\left<x\right>,-)$$ is seen to be the
colimit of representable functors.
\end{proof}

\begin{remark}
To summarize the above proof, the main ingredient is the observation
that the functor $\widetilde{pt}$ extends to $Super$ where it is
representable.  Furthermore, the representing object
$\zp\left<x\right>$ is a limit of a diagram of objects in $\Lambda$.
This representing diagram is not unique, however it does not prevent
us from easily proving Theorem \ref{functionsplus}, i.e. computing
the functions on $\widetilde{pt}$ below.

\end{remark}

\begin{cor}\label{pitaproof}
The functions on the $p$-adic superspace $\widetilde{pt}$ are
$\zp\left<x\right>$.
\end{cor}

\vspace{\baselineskip} \noindent\textbf{Acknowledgements.} We are
deeply indebted to V. Vologodsky for  his help with understanding of
the standard approach to the construction of Frobenius map and to M.
Kontsevich, N. Mazzari  and A. Ogus  for interesting discussions.
The second author would also like to thank the  Max Planck Institute
for Mathematics in Bonn where part of the work on this paper was
completed.

\noindent Department of Mathematics, University of California,
Davis, CA, USA \newline \emph{E-mail address}:
\textbf{schwarz@math.ucdavis.edu}

\noindent Department of Mathematics and Statistics, University of Windsor, Canada
\newline \emph{E-mail address}:
\textbf{ishapiro@uwindsor.ca}

\end{document}